\documentclass[12pt,psfig]{article}
\usepackage{amsthm,amssymb,amsmath,amsfonts,amscd}
\usepackage{graphicx,epsfig,latexsym}
\usepackage{cite,calc,xcolor,subfigmat,mathrsfs,dsfont}
\usepackage[
colorlinks,linkcolor=blue,citecolor=magenta]{hyperref}
\newcommand{\B}{\mathscr{B}}
\newcommand{\F}{\mathscr{F}}
\newcommand{\LST}{\mathscr{L}}
\newcommand{\E}{\mathscr{E}}
\newcommand{\Sl}{\mathfrak{sl}_2}
\newcommand{\T}{\mathfrak{T}}

\newcommand{\N}{\mathbb{N}}

\newcommand{\Span}{{\rm span}}
\newcommand{\ad}{{\rm ad}}

\newcommand{\0}{\mathbf{0}}
\newcommand{\C}{\mathcal{C}}

\newcommand{\IM}{{\rm im}\,}
\newcommand{\ie}{{\em i.e.,} }
\newcommand{\eg}{{\em e.g.,} }

\newcommand{\y}{\Delta}

\newcommand{\hot}{{\rm h. o. t. }}
\newcommand{\ddx}{\frac{\partial}{\partial x}}
\newcommand{\ddy}{\frac{\partial}{\partial y}}
\newcommand{\ddz}{\frac{\partial}{\partial z}}

\newtheorem{thm}{Theorem}[section]

\newtheorem{lem}[thm]{Lemma}

\theoremstyle{definition}
\newtheorem{defn}[thm]{Definition}
\newtheorem{exm}[thm]{Example}
\newtheorem{rem}[thm]{Remark}
\newtheorem{notation}[thm]{Notation}
\numberwithin{equation}{section}
\def\be {\begin{equation}}
\def\ee {\end{equation}}
\def\ba {\begin{eqnarray}}
\def\ea {\end{eqnarray}}
\def\bes {\begin{equation*}}
\def\ees {\end{equation*}}
\def\bas {\begin{eqnarray*}}
\def\eas {\end{eqnarray*}}
\def\bpr {\begin{proof}}
\def\epr {\end{proof}}

\begin{document}
\baselineskip=18pt
\renewcommand {\thefootnote}{\dag}
\renewcommand {\thefootnote}{\ddag}
\renewcommand {\thefootnote}{ }

\pagestyle{empty}

\begin{center}
\leftline{}
\vspace{-0.500 in}
{\Large \bf Normal forms of Hopf-zero singularity} \\ 

{\large Majid Gazor$^{*},$ Fahimeh Mokhtari\(\dagger\)}
\footnote{$^*\,$Corresponding author. Phone: (98-311) 3913634; Fax: (98-311) 3912602; Email: mgazor@cc.iut.ac.ir; \(\dagger\,\) Email: f.mokhtari@math.iut.ac.ir. }

\vspace{0.105in} {\small {\em Department of Mathematical Sciences,
Isfahan University of Technology
\\[-0.5ex]
Isfahan 84156-83111, Iran }}

\today

\vspace{0.05in}

\noindent
\end{center}

\vspace{-0.10in}

\baselineskip=15pt

\:\:\:\:\ \ \rule{5.88in}{0.012in}

\begin{abstract}

The Lie algebra generated by Hopf-zero classical normal forms is decomposed into two versal Lie subalgebras.
Some dynamical properties for each subalgebra are described; one is the set of all volume-preserving conservative systems while the other is the maximal Lie algebra of nonconservative systems. This introduces a unique conservative--nonconservative decomposition for the normal form systems.
There exists a Lie--subalgebra that is Lie-isomorphic to a large family of vector fields with Bogdanov--Takens singularity. This gives rise to a conclusion that the local dynamics of formal Hopf-zero singularities is well-understood by the study of Bogdanov--Takens singularities. Despite this, the normal form computation of Bogdanov-Takens and Hopf-zero singularities are independent. Thus, by assuming a quadratic non-zero condition, complete results on the simplest Hopf-zero normal forms are obtained in terms of the conservative--nonconservative decomposition. Some practical formulas are derived and the results implemented using Maple.
The method has been applied on the R\"{o}ssler and Kuramoto--Sivashinsky equations to demonstrate the applicability of our results.

\vspace{0.10in} \noindent {\it Keywords:} \ Normal form;
Hopf-zero singularity; \(\Sl\)-representation; Conservative and nonconservative decomposition.

\vspace{0.10in} \noindent {\it 2010 Mathematics Subject Classification}:\, 34C20; 34A34.

\noindent \rule{5.88in}{0.012in}
\end{abstract}

\vspace{0.2in}

\section{Introduction}
In this paper we are concerned with computing the simplest normal form of the system
\be\label{Eq1} \dot{x}:= \hot, \quad \dot{y}:= z+ \hot, \quad \dot{z}:= -y+ \hot, 
\ee where \(\hot\) denotes formal nonlinear terms (higher order terms) with respect to \((x,y,z)\in \mathbb{R}^3.\) Our normal form computation is a local tool and can be used for local dynamics analysis of formal flows; this is because the convergence of transformations are not discussed here. The system (\ref{Eq1}) can
be transformed into the first level (classical) normal form
\ba\label{CNF}
\left\{
  \begin{array}{ll}
 \frac{dx}{dt}=
\sum^\infty_{i+2j=2} a_{ij}
x^{i}{(y^2+z^2)}^{j}, &
\\
\frac{dy}{dt}= z
 + \sum^\infty_{i+2j=1} x^i{(y^2+z^2)}^{j}(b_{ij}y+ c_{ij}z\big),&
\\
\frac{dz}{dt}= -y+ \sum^\infty_{i+2j=1} x^i{(y^2+z^2)}^{j}(b_{ij}z- c_{ij}y\big),&
  \end{array}
\right.
\ea
where \(a_{ij}, b_{ij},c_{ij}\in \mathbb{R}\); see also \cite[Equation 2.15]{AlgabaHopfZ}. All the existing results on the simplest normal forms of this singularity have only dealt with the cases where \(a_{20}b_{10}\neq 0\); see \cite{AlgabaHopfZ,ChenHopfZ03,ChenHopfZ,YuHopfZero,Mokhtari}. The main reason for this is that terms corresponding to \(a_{01},\) \(b_{10},\) \(c_{10},\) and \(a_{20}\) all have the same grades in the usual graded structures; \ie the gradings may not distinguish these four monomial vector fields. Therefore, one needs to deal with all four terms in the computations which is a tremendously difficult job. A novelty in the results of Algaba {\em et. al.} \cite{AlgabaHopfZ} was to notice that the first two non-zero terms (associated with \(a_{20}\) and \(b_{10}\)) play the key role in the computations and thus, generic conditions with respect to \(a_{20}\) and \(b_{10}\) are assumed. Then, they obtained the simplest normal forms and orbital normal forms. Chen {\em et. al.} \cite{ChenHopfZ03,ChenHopfZ} approached this family of systems using an essentially different method and provided an independent solid proof for normal form uniqueness. Yu and Yuan \cite{YuHopfZero} made an efficient computer program to compute the simplest normal form of this family of systems.

Some discussions on convergence and divergence of normal forms have been made in \cite{GazorMokhtariInt}. For instance, we computed the numerically suggested radius of convergence of the second level normal forms associated with volume-preserving Hopf-zero singularity. We do not address the convergence problem in this paper; see \cite{IoosLombardi,Stroyzyna,StZolDiv,Zoladek02} for some related results. Therefore, all our claims with regards to dynamics analysis are limited to the formal flows of formal systems.

\pagestyle{myheadings} \markright{{\footnotesize {\it M. Gazor and F. Mokhtari \hspace{3.2in} {\it Hopf-zero singularity}}}}

A motivation of this paper is to use the monomial term corresponding to \(a_{01}\) as the main player in calculations; this is to complete the existing results on this problem with a quadratic non-zero term. (Throughout this paper the only assumption is \(a_{01}\neq 0.\)) The grading in our approach distinguishes this vector field (grade \(1\)) from the other three (grade \(2\)). This is achieved through a \(\Sl\)-representation for the classical normal forms. This greatly simplifies the computations; we obtain complete results on the simplest normal form of these systems without any extra generic conditions. Our approach can be applied on many well known models such as R\(\ddot{\hbox{o}}\)ssler and Kuramoto--Sivashinsky equations.

The computational burden has been the main obstacles of most classification problems in the normal form literature. A systematic approach is required for such computations; a Lie-graded structure is an important instrument where grade-homogenous parts of vector fields are simplified inductively. In this process, a basis for the space of grade-homogenous vector fields is applied. The mostly used choice for the basis has been the monomial vector fields oblivious of their dynamics. This usually involves large matrices and it is hard to find the patterns of computations. The conservative and nonconservative polynomial vector fields play this role in this paper instead of the monomials. Then, computations yield more solid patterns compared to when monomial vector fields are used. Indeed both transformation generators and normal forms are presented through conservative and nonconservative grade-homogeneous vector fields. This is a new feature of this paper that distinguishes our results from the existing results on the simplest normal forms of this singularity.
This is accomplished via a \(\Sl\)-representation for the classical normal form vector fields. This technique has been mainly applied on nilpotent singularities. Thus, it may be surprising to see that the theory of \(\Sl\)-representation is applied to a non-nilpotent system. Another novelty of our results is to use \(\Sl\)-style for the second level normal forms; \(\Sl\)-styles have only been used for the first level normal forms in the existing literature; see also \cite{BarrioPalacian,PalacChaos}.

Any Lie (sub)algebra structure may have interesting dynamics interpretations and normal form theory provides a powerful tool for such descriptions.
For our first instance, recall that a Lie subalgebra generates a group of transformations and the group establishes an equivalence relation. Hence, it gives us a classification within the space through infinite level normal forms. These usually introduce important families of vector fields.
Here, we denote \(\LST\) for the Lie algebra generated by the first level Hopf-zero normal forms. Two transversal Lie subalgebras for \(\LST\) are presented. These represent the quasi-Eulerian vector fields and volume-preserving vector fields with a first integral. The second and more interesting example is described as follows. There exists a subalgebra from \(\LST\) that is Lie-isomorphic to a subalgebra \(\LST_{b}\) from the Lie algebra generated by all two dimensional vector fields with Bogdanov--Takens singularity; see Theorem \ref{23e}. This gives rise to the fact that local dynamics of any planar reduced (by ignoring the phase coordinate) system from Hopf-zero singularity can be embedded into the flow of a Bogdanov-Takens singularity; see Theorem \ref{BTintoHopf0}. Therefore, the local planar flow associated with Hopf-zero is well-understood by studying that of Bogdanov-Takens; see our further detailed discussion following Remark \ref{Rem26}. In other words, the local reduced system of Hopf-zero holds less complexity than Bogdanov-Takens.
The Lie isomorphism further provides an explanation for why the \(\Sl\)-representation works fine for non-nilpotent singularities and also suggests that our techniques may be applicable to some other non-nilpotent singularities.

The rest of this paper is organized as follows. The conservative-nonconservative decomposition for normal forms are introduced in Section \ref{secDec}. This is achieved by a \(\Sl\)-representation for \(\LST\) and presenting two transversal Lie subalgebras. Some properties for each family are described. The second level normal form is computed in Section \ref{2NF}. In Section \ref{SNF} we obtain the simplest normal forms. The procedure is divided into three cases and accordingly their simplest normal forms are computed in three subsections. The results are applied on R\"{o}ssler and Kuramoto--Sivashinsky equations in Section \ref{sec4}.

\section{Lie algebra \(\LST\) and its \(\Sl\)-representation } \label{secDec}

The space of all classical normal forms governed by Equation (\ref{CNF}) is denoted by \(\LST\). This section provides a \(\Sl\)-representation for \(\LST\). Any column vector \([f_1, f_2, f_3]^T\) is associated with the vector field \(v:=f_1\ddx+ f_2\ddy+f_3\ddz\) and vice versa. Hence, \(v\) generates a system given by \([\dot{x}, \dot{y}, \dot{z}]^T:=[f_1, f_2, f_3]^T.\) Besides, \(v\) acts as a differential operator on scalar functions, say \(g(x, y, z)\), defined by
\bes
v(g):= \frac{f_1\partial g}{\partial x}+ \frac{f_2\partial g}{\partial y}+\frac{f_3\partial g}{\partial z}.
\ees Thereby, terminologies of ``system'', ``vector field'' and ``differential operator'' are interchangeably used. For any vector field \(w,\) we define
\bes
wv:= w(f_1)\ddx+ w(f_2)\ddy+w(f_3)\ddz.\ees
Now \(\LST\) is a Lie algebra by \(\ad_v(w):=[v,w]= vw-wv\) for any \(v,w\in \LST.\)
Any \(\Sl\) Lie algebra is represented by a triad \(\{N, M, H\}\). In this paper they are introduced by
\bas
N&:=& (y^2+z^2)\ddx,\\
M&:=&-\frac{xy}{(y^2+z^2)}\ddy- \frac{xz}{(y^2+z^2)}\ddz,\\
H&:=&2x\ddx-y\ddy- z\ddz.
\eas It is easy to verify \({[M, N]}= H,\) \({[H, M] }= 2M,\) and \({[H, N]}= -2N.\) The triad differential operators \(\{\ad_M, \ad_N, \ad_H\}\) act on the space of vector fields \(\LST.\) However, \(\LST\) is not invariant under this action. Any vector field from \(\LST\) generates an orbit under this action and the orbit is terminated when it reaches either the zero vector field or a non-permissible vector field (not defined in \(\LST\)).
The idea is to construct certain orbits from which the vector fields from \(\LST\) are decomposed in terms of the orbit elements.
This is done by the following definition. This provides a tool to use a similar approach to the method applied by Baider and Sanders  \cite{baiderchurch,baidersanders,BaidSand91}.

\begin{defn}
Define {\bes F^{-1}_0:= 2{(y^2+z^2)}\ddx, \, E^0_0=x\ddx+ \frac{1}{2}y\ddy+\frac{1}{2} z\ddz, \Theta^0_0:= z\ddy- y\ddz,\ees} and
\ba\label{orbits1}
{F}^{l}_k&:=&\frac{(-1)^{l+1}(k-l+1)!}{2^{l+1}(k+2)!}\ad_M^{l+1}(y^2+z^2)^kF^{-1}_0, \qquad\qquad\qquad \hbox{ for } -1\leq l\leq k,
\\\label{orbits2}
{E}^l_k&:=&\frac{(-1)^l(k-l)!}{2^lk!} \ad_M^l(y^2+z^2)^kE^0_0, \qquad\qquad\quad\qquad \qquad\ \ \hbox{ for } 0\leq l\leq k,
\\\label{orbits3}
\Theta^l_k&:=&\frac{(-1)^l(k-l)!}{2^lk!} \ad_M^l(y^2+z^2)^k\Theta^0_0, \qquad\qquad\qquad\qquad\quad\ \ \hbox{ for } 0\leq l\leq k,
\ea where \(\ad_M v:= [M, v]\) and \(\ad_M^nv:= [M, \ad_M^{n-1} v]\) for any natural number \(n\).
\end{defn}
These give rise to the following theorem.
\begin{lem}\label{ThmLie} The formulas for \(F^l_k\), \(E^l_k\) and \(\Theta^l_k\) in terms of \((x, y, z)\)-coordinates are given by
\ba\nonumber
F^l_k&=& x^{l} {(y^2+z^2)}^{k-l}\left((k-l+1) x\ddx-
\frac{(l+1)}{2} y\ddy-\frac{(l+1)}{2} z\ddz\right),
\\\label{Eul}
E^l_k&=& x^{l}{(y^2+z^2)}^{k-l}\left(x\ddx+ \frac{1}{2}y\ddy+\frac{1}{2}z\ddz\right),
\\\nonumber
\Theta^l_k&=& x^{l}(y^2+z^2)^{k-l}\left(z\ddy-y\ddz\right).
\ea
\end{lem}
\bpr
The proof readily follows an induction on \(l.\)
\epr
\noindent A straightforward calculation proves the following lemma.

\begin{lem}\label{structure} The structure constants for the Lie algebra \(\LST\) is governed by
\bas\label{aa}
{[F^l_k, F^{m}_{n}]}&=&
\big((m+1)(k+2)-(l+1)(n+2)\big)F^{l+m}_{k+n},\\
{[F^l_k, E^{m}_{n}]}&=&
{\frac{(n+2)\big(m(k+2)-n(l+1)\big)}{(k+n+2)}}E^{l+m}_{k+n}-
{\frac{k(k+2)}{k+n+2}}F^{l+m}_{k+n},\\
{[F^l_k,\Theta^m_n]}&=&\big(m(k+2)-n(l+1)\big)\Theta^{l+m}_{k+n},
\\
{[E^l_k, E^{m}_{n}]}&=& (n-k)E^{l+m}_{k+n},
\\
{[E_k^l, \Theta^m_n]}&=&n\,\Theta^{l+m}_{ k+n},
\\
{[\Theta_k^l,\Theta_n^m]}&=&0.
\eas
\end{lem}

Let
\bes
\F:=\Span\left\{a_0F^{-1}_0+\sum a^l_k F^l_k\,|\, 1\leq k, -1\leq l\leq k, a^l_k\in \mathbb{R}\right\},
\ees  and
\bes
\T:= \Span\left\{\Theta^0_0+\sum c^l_k \Theta^l_k\,|\, 1\leq k, 0\leq l\leq k, c^l_k\in \mathbb{R}\right\}.
\ees Hence, any \(v\in \T\) has only phase components in cylindrical coordinates. The set of all formal first integrals for \(v\in \T\) is the algebra generated by \(x\) and \(y^2+z^2\); see \cite[Lemma 2.1]{GazorMokhtariInt}. Further, denote
\be\label{Lie}
\E:=\Big\{\sum b^l_{k} E^l_k \;|\;  {l+k\geq 1}, k\geq l,
b^l_{k}\in \mathbb{R}\Big\}. \ee

The following theorem indicates that the space of all volume-preserving and conservative normal forms is given by \(\F\oplus\T\); see \cite{Wiggins3D} for relevant results on three dimensional volume-preserving vector fields and their normal forms. Furthermore, it recalls that nonzero vector fields from \(\E\) are nonconservative.

\begin{thm}\label{itemsThm} The following holds.
\begin{enumerate}
\item\label{111} Any differential system governed by Equation \eqref{CNF} is associated with a vector field \(v\in \LST\) and vice versa.
\item \label{23f} The vector spaces \(\F, \E,\) and \(\T\) are transversal Lie subalgebras in \(\LST,\) \ie \bes\LST=\F\oplus \E\oplus \T.\ees
\item \label{23d} The Lie subalgebra \(\T\) is a Lie ideal for \(\LST.\) Thus, \(\F\oplus\T\) is also a Lie subalgebra for \(\LST.\)
\item \label{231a} The algebra of first integrals for any nonzero element from \(\T\) is \(<x, y^2+z^2>.\)
\item \label{23a} Let
\bes v= \Theta^0_0+a_0 F^{-1}_0+ \sum\left( a^l_k F^l_k+ c^l_k \Theta^l_k\right)\in \LST, \ \hbox{ where } a_0\neq 0.\ees
Then, there exists a unique formal first integral \(f=\sum  a^l_k f^l_k\) (modulo scalar multiplications) such that the algebra of first integrals for
\(v\) is \(\langle f\rangle,\) where \(f^l_k:= x^{l+1}(y^2+z^2)^{k-l+1}.\)
\item \label{23b} There is no nonzero first integral for any nonzero vector field from \(\E.\) Furthermore, nonzero vector fields from \(\E\) are not volume preserving.
\item \label{23c} The space \(\F\oplus\T\) is the maximal vector space of volume-preserving vector fields. The algebra of first integrals for any \(v\in \F\oplus\T\) is nontrivial.
\item \label{23sdynm} The flow generated by any \(v\in \T\) is static in the amplitude and \(x\)-coordinates. In other words, it can only be dynamic in the
phase coordinate.
\item \label{23static} For any \(v\in \F\oplus \E\), the generated flow is static in the phase coordinate.
\end{enumerate}
\end{thm}
\bpr The space \(\LST\) is defined such that claim \ref{111} holds.
The proofs for \ref{23f} and \ref{23d} follow Lemma \ref{structure}.
The claims \ref{231a}, \ref{23a} and \ref{23b} are proved by
\cite[Lemma 2.1]{GazorMokhtariInt}, \cite[Proposition 2.2]{GazorMokhtariInt} and \cite[Theorem 2.3]{GazorMokhtariEul}, respectively.
Since \(\LST= \E\oplus \F\oplus\T,\) the proof of \ref{23c} is straightforward. Claim \ref{23sdynm} is true because \(\T\) is generated by all vector fields from \(\LST\) with zero \(x\) and amplitude components.
Since \(\F\oplus \E\) is transversal to the Lie algebra \(\T\), the flow generated by \(v\in \F\oplus \E\) is static in phase coordinate.
\epr

\begin{thm}\label{23e}
The space \(\F\oplus\E\) is a Lie algebra and there exists a Lie-isomorphism \(\psi\) to a proper Lie subalgebra \(\LST_b\) of the Lie algebra generated by two dimensional Bogdanov--Takens singularities.
\end{thm}
\bpr
Denote \(\LST_{\rm BT}\) for the Lie algebra generated by all vector fields of the form
\be\label{BTEq}
w:= -\overline{x}\frac{\partial}{\partial \overline{y}}+\hot.
\ee The negative sign is chosen such that it matches with the notation of \cite{BaidSand91}. Define a map
\ba\label{h}
&\psi: \F\oplus\E \rightarrow \LST_{\rm BT},&
\ea governed by
\bes
F^l_k\mapsto (k+2)A^{k-l}_k,
\
E^l_k\mapsto B^{k-l}_k,
\ees where \(A^l_k\) and \(B^l_k\) are defined by Baider and Sanders \cite[Equations 3.6a and 3.6b]{BaidSand91}.
 Comparing the two sets of structure constants, \(\psi\) is a Lie-monomorphism and the claim follows by \(\LST_b:= \psi(\F\oplus\E)\).
\epr
\begin{rem} \label{Rem26}
The claim \ref{23a} from Theorem \ref{itemsThm} implies that \(f^l_k(x, y, z):= x^{l+1} (y^2+z^2)^{k-l+1}\) refers to \(F^l_k\) in \(\LST\) and from \cite{BaidSand91} the monomial \(h^{k-l}_k(\overline{x}, \overline{y}):= \overline{x}^{k-l+1}\overline{y}^{l+1}\) is associated with \(A^{k-l}_k\) in \(\LST_b.\) The monomial \(f^{k+1}_k=x^{k+2}\) does not produce a permissible vector field in \(\LST\). However, \(h^{-1}_k=\overline{y}^{k+2}\) leads to \(A^{-1}_k.\) Here the \(\Sl\)-orbits given by \eqref{orbits1} are terminated when they lead to non-permissible vector fields (\ie \(l\) is increased to \(k+1\)), while the \(\Sl\)-orbits given by \cite[Equations 3.6a and 3.6b]{BaidSand91} are terminated by zero (\ie \(A^{-2}_k=0\)).
\end{rem}

Now we conclude that the local dynamics of a Hopf-zero normal form is well-understood by the study of Bogdanov--Takens singularities. The basis of our claim is as follows. There is a common practice to ignore the phase component of the Hopf-zero normal form and obtain a planar reduced system. Then, study the dynamics of the planar reduced system and extract the full three dimensional dynamics from the planar reduced system; see \cite{DumortierHopfZ,AlgabaHopfElect,HarlimLangf,LangfTori,LangfHopfSteady,LangfHopfHyst}. Given this, we prove that for any planar reduced system obtained from a Hopf-zero normal form, there exists a Bogdanov--Takens singularity from \(\LST_{b}\) such that the flow of the planar reduced system is embedded into the flow of Bogdanov--Takens singularity; see Theorem \ref{BTintoHopf0}. It is interesting to note that the dynamics of Bogdanov--Takens singularities are expected to be more rich than that of the planar reduced systems obtained from Hopf-zero. This is because of two reasons. First, the embedding between the two flows assigns the square of the amplitude variable from the reduced systems into a state variable from Bogdanov--Takens systems that is not necessarily nonnegative. This excludes certain dynamical complexities. Second (and more important than the first), a complement space to \(\LST_b\) is expected to absorb most of the dynamics of Bogdanov--Takens singularity. This is because all nonzero nonlinear terms from \(\LST_b\) appearing in a Bogdanov--Takens singularity are simplified in the \(\Sl\)-classical normal forms; see Remark \ref{Rem32}. This contributes to the possible exclusion of many complex dynamical behaviors.

\begin{thm}\label{BTintoHopf0}
For any Hopf-zero classical normal form, there exists a Bogdanov--Takens singularity such that the flow generated by the planar reduced system (associated with Hopf-zero) is embedded into the Bogdanov--Takens' flow.
\end{thm}
\bpr Consider the vector fields from \(\LST\) in cylindrical coordinates. Now for any \(v\in \LST,\) let
\bes v:= \tilde{v}+\hat{v}, \hbox{ where } \tilde{v}\in \T \hbox{ and } \hat{v}\in \F\oplus\E.\ees
Define \(w:= \psi(\hat{v})\in \LST_b\) and \(\check{v}\) as the two dimensional vector field obtained from \(\hat{v}\).
Denote nonnegative real numbers by \(\mathbb{R}^+.\) The changes of variables
\ba\nonumber
&\check{\psi}_{\theta_0}: \mathbb{R}\times \mathbb{R}^+\rightarrow\mathbb{R}^+\times \mathbb{R}\subset \mathbb{R}^2, &
\\\label{hhomeo}
&\check{\psi}_{\theta_0}(x, \rho)=(\overline{x}, \overline{y}), \ \overline{x}(x, \rho):= \rho^2 \hbox{ and } \overline{y}(x, \rho):=x, &
\ea is a homeomorphism transforming \(\check{v}\) into \(w\). Hence, \(\check{\psi}\) embeds the flow of \(\check{v}\) into the flow of \(w\). This completes the proof.
\epr

\begin{notation}
Throughout this paper we use Pochhammer \(j\)-symbol notation, that is,
\be\left(a\right)^j_b:=\prod_{i=0}^{j-1}(a+ib), \ee
for any natural number \(j\) and real number \(b.\) Further, for integer numbers \(p, q, r,\) we denote \be p\equiv_{_{q}} r\ee when there exists an integer \(k\) such that \(p-r= kq.\)
\end{notation}

Assume that \(v=v_0+v_1+v_2+\cdots,\) where \(v_k\in \LST_k\) and \(\LST_k\) denotes homogenous vector fields with grade \(k\) for any \(k.\) Let
\bas
&&L^{n,1}:\LST_n\rightarrow \LST_n, \\
&&\;\; S_n\mapsto [S_n, v_0].
\eas Next we inductively define
\bas
&&L^{n,N}:\LST_n\times \ker L^{n-1,N-1}\rightarrow \LST_n,\\
&&(Y_{n}, Y_{n-1}, \ldots , Y_{n-N+1})\mapsto \sum^{N-1}_{i=0} [Y_{n-i}, v_{i}].
\eas There exists a complement space \(\C^{n,N}\) such that \(\IM L^{n,N}\oplus \C^{n,N}=\LST_n.\) The vector field
\bes w=w_0+w_1+w_2+\cdots\ees
is called an \(N\)-th (infinite) level normal form when \(w_n\in \mathcal{C}^{n,N}\) (\(w_n\in \mathcal{C}^{n,n}\)) for any \(n.\)
Then, there exist invertible transformations such that \(v\) can be transformed into the \(n\)-th and infinite level normal form; see \cite{GazorYuSpec}.

\section{Hypernormalization } \label{2NF}

In this section we provide a hypernormalization for the system (\ref{Eq1}). Any such system is transformed into Equation \ref{CNF} and then by the item \ref{111} of Theorem \ref{itemsThm}, it can be represented by
\be\label{eq28} v^{(1)}= \Theta^0_0+a_0 F^{-1}_0+ \sum\left( a^l_k F^l_k+ b^l_k E^l_k+ c^l_k \Theta^l_k\right)\in \LST.\ee

\begin{lem}\label{2nd} For any vector field \(v^{(1)}\) given by equation \eqref{eq28} where \(a_0\neq 0,\)
there exist invertible changes of state variables transforming \(v^{(1)}\) into
\be\label{v1} v^{(2)}=
\Theta^0_0+a_0F^{-1}_0+\sum^\infty_{k=1}\left(a_k{F^k_k}+b_k{E^k_k}+c_k{\Theta^k_k}\right).
\ee
\end{lem}
\bpr
Define the grading function by \(\delta(F^l_k)=\delta(E^l_k)=\delta(\Theta^l_k)=k\), then proof readily follows \([E^l_k,F^{-1}_0]= -2lE^{l-1}_k,\) \({[F^l_{k},F^{-1}_0]}=-2(l+1)F^{l-1}_{k}\) and \([\Theta^l_k,F^{-1}_0]=-2l\Theta^{l-1}_k\).
\epr

A bounded and periodic transformation \(\varphi\) is given by
\be\label{LinChange}[x(t), \rho(t), \Theta(t)]=\varphi\big(x(t), \rho(t), \theta(t)\big):=[x, \rho, \theta-t],\ee
where \(\Theta\) and \(\theta\) are denoted for the new and old phase variables, respectively. The change of variable in \(\varphi\) omits \(\Theta^0_0\) from
the vector field \(v^{(2)}\). Once the simplest normal form computation is exhausted, the map \(\varphi^{-1}\) adds \(\Theta^0_0\) back into the system; also
see \cite[Theorem 4.1]{GazorMokhtariEul}. This technique has been used in perturbation theory; see \cite[Lemma 5.3.6]{MurdBook}. Since the linear part of the system \(\Theta^0_0\) commutes with the classical normal form, its elimination (using \eqref{LinChange}) does not change the formal normal form. However, if one uses a truncation at a grade (say \(k\)) plus a remainder, then the transformed system is nonautonomous beyond grade \(k.\) Another common way of looking at normal forms is to consider a formal normal form as smooth modulo a flat vector field; \eg see \cite{Murd08}. This gives rise to an autonomous formal normal form plus a non-autonomous flat system. This was brought to our attention by James Murdock. This is new in the normal form literature and substantially reduces the computational burden.

After eliminating \(\Theta^0_0,\) we may change \(a_0\) to any arbitrary number \(\tilde{a}_0\) via the time rescaling \(t:=\frac{\tilde{a}_0}{a_0}\tau\) where \(t\) and \(\tau\) denote the old and new time variables, respectively. Therefore, without the loss of generality we may assume that
\be\label{v2} v^{(2)}= F^{-1}_0+\sum^\infty_{k=1}\left(a_k{F^k_k}+b_k{E^k_k}+c_k{\Theta^k_k}\right).
\ee

\begin{rem}\label{Rem32}
All nonlinear Hamiltonian and Eulerian terms that appear in the first level normal form of Bogdanov--Takens singularity in \cite{BaidSand91} are of the form \(A^{-1}_k\) and \(B^0_k\) while these terms do not belong to \(\LST_b\); see Theorem \ref{23e}. Therefore, the computations in hypernormalization steps in this paper are essentially independent from \cite{BaidSand91} despite similarity of the procedures and the existence of a Lie--isomorphism between \(\LST_b\) and \(\F\oplus \E\).
\end{rem}

For hypernormalization we need to use a normal form style. A style is a rule on how to choose complement spaces in the normal form computation. We use the \(\Sl\)-style in a hypernormalization step. The \(\Sl\)-style in \cite{BaidSand91} states that the only nonlinear terms from \(\ker \ad_{A^{-1}_0}\) can stay in the classical normal forms. Those are \(A^{-1}_k\) and \(B^0_k\) for any \(k\); see \cite{BaidSand91}. Lemma \ref{2nd} proves that the terms \(F^k_k, E^k_k,\) and \(\Theta^k_k\) (see Equations \eqref{orbits1}--\eqref{orbits3}) may stay in our second level normal form system. These are described by \(\ker \ad_{F^1_0},\) if we would define \(F^1_0.\) Indeed, we embed \(\LST\) into a bigger algebra (\ie as a proper Lie subalgebra), where the bigger algebra contains \(F^1_0\). As far as algebra is concerned, this is simply performed via a simple generalization of Lemma \ref{structure}. However, this extension may not have a justification in nonlinear dynamics; \ie assigning any term in an ODE system or a near-identity transformation to those extra generated vector fields (beyond \(\LST\)) is not permissible. Indeed, from the bigger algebra, \(F^1_0\) is only used in definitions of \(\mathscr{Gamma}\) given by Equations \eqref{gamma}, \eqref{gamma2} and \eqref{gamma3}. The map \(\mathscr{Gamma}\) is merely an instrument that substantially facilitates the computation. The original ideas behind definition of \(\mathscr{Gamma}\) and our gradings come from \cite{BaidSand91}.

Once the second level normal form is calculated, the \(\Sl\)-style is extended to a formal basis style. Formal basis style uses an order on normal form terms of \eqref{v1} to distinguish between alternative terms for elimination. Here, we give priority of elimination to conservative terms over nonconservative terms of the same grade in style \(I\) and vice versa in style II; see \cite{GazorYu,GazorYuSpec,GazorYuFormal,Murd04,Murd08} for further information on formal basis style. Let
\bes \mathscr{B}_N=\left\{F^{l}_{k}, \,E^m_n,\Theta^m_n|\,-1\leq {l}\leq{ k},\,0\leq m\leq n,\delta({F^{l}_{k}})=\delta(E^m_n)=\delta(\Theta^m_n)=N\right\}\ees be a basis for the vector space \(\LST_N\) and \(\delta\) denote for a grading function.
Two formal basis styles are defined through the following orderings on \(\B_N\) and are used in this paper.
\begin{itemize}
 \item Style \({I}\): \(E^l_k\prec F^m_n\prec\Theta^i_j\) for any \(l, k, m, n, i, j\). Furthermore, \(F^l_k\prec F^m_n,\)
\(E^l_k\prec E^m_n\) and \(\Theta^l_k\prec\Theta^m_n\) if \(k<n.\) This gives a priority for elimination of \(E^l_k\)-terms rather than \(F^m_n\)-terms.
 \item Style \({II}\): \( F^m_n\prec E^l_k\prec\Theta^i_j\) for any \(l, k, m, n, i, j\). Besides, \(F^l_k\prec F^m_n,\)
\(E^l_k\prec E^m_n\) and \(\Theta^l_k\prec\Theta^m_n\) if \(k<n.\) The priority of elimination is with \(F^m_n\)-terms rather than \(E^l_k\)-terms.
\end{itemize}

\section{The simplest normal forms}\label{SNF}

Throughout this paper we assume that there exist \(a_l\neq0\) and \(b_k\neq0\) for some \(l\) and \(k.\)
Define
\be {r}:= \min\{i\,|\,a_i\neq 0, i\geq1\},  s:= \min\{j\,|\,b_j\neq 0,j\geq1\}, \hbox{ and } p:= \min\{j\,|\,c_j\neq 0,j\geq1\}.\ee
Then, in order to obtain complete results for the simplest normal forms we divide the problem into the following three cases:
\be\label{3cases} (i)\ \ {r}<s, \ \ \ \quad (ii)\ \ {r}>s, \  \ \ \quad (iii)\ \ {r}=s.
\ee

In this section we compute the simplest normal forms for the three cases \eqref{3cases} in the following three subsections. We hereby acknowledge N. Sadri's help for an independent verification, and detecting a few errors, of the formulas.

\subsection{Case i: \({r}<s\).}\label{casei}

Assume that
$$v^{(2)}:=F^{-1}_0+{a_{r}}F^{r}_{r}+\sum_{k={r}+1}^{\infty}a_k{F^k_k}
+\sum^\infty_{k=s}b_k{E^k_k}+\sum^\infty_{k=1}c_k{\Theta^k_k},$$ and define \be\label{grad11}\delta(F^l_k)= \delta(E^l_k)={r}(k-l)+k \hbox{  and } \delta(\Theta^l_k)= {r}(k-l+1)+k+1.\ee

The vector field
\be\label{Fr}
\mathbb{F}_{r}:=F^{-1}_0+{a_{r}}F^{r}_{r}
\ee plays an important role in further normalization of $v^{(2)}.$
For any arbitrary \(\alpha\in \mathbb{R}, \alpha>0, \) through changes of variables
\bes t:=\left(\frac{\alpha\, {\rm sign}(a_r)}{a_r}\right)^{\frac{1}{r+1}}\tau, \
x:= \left(\frac{\alpha\, {\rm sign}(a_r)}{a_r}\right)^{\frac{1}{r+1}}X, \ y:= Y, \ \hbox{ and } \ z:= Z,\ees  we can change \(a_r\) into
\(\alpha \,{\rm sign}(a_r).\) Thereby, without the loss of generality we assume that \(a_r\) in Equation \eqref{Fr} is a non-algebraic number. When it comes to practical normal form computation by using a computer, this assumption is not valid; computers do not recognize irrational numbers.

\begin{lem}
The \(({r}+1)\)-th level normal form of (\ref{Eq1}) is
\be v^{({r}+1)}: =
F^{-1}_0+{a_{r}}{F^{r}_{r}}+{b_s}^{({r}+1)}{E^s_s}+\sum_{k>{r}} a_k^{({r}+1)}{F^k_k}+
\sum_{k>s} b_k^{{({r}+1)}}{E^k_k}+\sum_{k\geq1}{c}_k^{({r}+1)} \Theta^k_k, \ee where in
\begin{itemize}
  \item style I, \( a_k=0\) for \({k\equiv_{ 2({r}+1)} {2r}},\) and \({k\equiv_{ 2({r}+1)} {{r}}}\), and  \(b_k=0\) for \({k\equiv_{ 2({r}+1)} {-1}}\), while \( c_k=0\) for \({k\equiv_{ 2({r}+1)} {-1}}.\)
  \item style II, \( a_k=0\) for \({k\equiv_{ 2({r}+1)} {2r}}\), \({k\equiv_{ 2({r}+1)} {{r}}},\) and \({k\equiv_{ 2({r}+1)} {-1}}\), where \( c_k=0\) for \({k\equiv_{ 2({r}+1)} {-1}}.\)
\end{itemize}
\end{lem}
\bpr
It is easy to see that \(v^{(2)}\in {\ker\ad_{F^{1}_0}}.\) Thereby, in order to compute the higher level normal forms we follow \cite{GazorMoazeni,baidersanders,GazorMokhtariEul,GazorMokhtariInt} and define
\be\label{gamma}\mathscr{Gamma}:= \ad_{F^{1}_{0}}\circ \ad_{\mathbb{F}_{r}}.\ee Then,
\bas
\mathscr{Gamma} (F^{l}_{k})&=&-4(l+1)(k-l+2)F^{l}_{k}+2a_{{r}}(k-l+1)
\big((k-l)({r}+1)+r-l\big) F^{{r}+l+1}_{{r}+k},
\\
\mathscr{Gamma}(E^l_k)&=&-4l(k-l+1)E^l_k
+\frac{2a_{{r}}(k+2)(k-l)\big((k-l)({r}+1)-l\big)}{{r}+k+2}E^{{r}+l+1}_{{r}+k}\\
&&+\frac{2a_{{r}}{r}({r}+2)(k-l+1 )}{{r}+k+2}F^{{r}+l+1}_{{r}+k},
\\
\mathscr{Gamma}(\Theta^{l}_{k})&=& -4l(k-l+1)\Theta^{l}_{k}+2a_{{r}}(l-k)\big(l({r}+2)-k({r}+1)\big)
\Theta^{{r}+l+1}_{{r}+k}.
\eas
Given the ordering for \(\B_N,\) the matrix representation of \(\mathscr{G}\) is lower triangular. Thus,
\be
\ker(\mathscr{Gamma})= \Span \left\{\mathcal{F}^{-1}_{k,r},\mathcal{E}^{0}_{k,r}, {\mathcal{X}_{r}^{k}}, \mathcal{T}^0_{k,r}\,|\, k\in \N\right\},
\ee
where \(\mathcal{F}^{-1}_{k,r} , \mathcal{X}_{r}^{k}, \mathcal{T}^0_{k,r}\) are defined in \cite[Equation 3.4]{GazorMokhtariInt}. In addition
\ba
\mathcal{E}^{0}_{k,r}&:=&\sum_{m=0}^{k}\frac{{a_{{{r}}}}^{m}(k+2)(k)^m_{-2}}{m!2^m(m{r}+k+2)}E^{m({r}+1)}_{m{r}+k}
-\sum_{m=0}^{k}{a_{{{r}}}}^{m}{e}_{{k,}m}F^{m({r}+1)}_{m{r}+k},
\ea here \(\delta(\mathcal{E}^{0}_{k,r})\equiv_{({r}+1)}0,\) \( e_{{k,}0}=0,\) and
 \bas
e_{{k,}m+1}&=&{\frac{\big({r}+(k-2m)({r}+1)\big)e_{k,m}}{2\big((m+1)({r}+1)+1\big)}}-
{\frac{{r}({r}+2)(k+2)(k)_{-2}^{m}}{m!2^{m+1}\big((m+1)({r}+1)+1\big){\big(m{r}+k+2\big)^2_r}}.
}\eas
 On the other hand
 \ba\nonumber
\\\label{k3}
{[\mathcal{E}^{0}_{2k,r}, {\mathbb{F}_{r}}]}&=&{{a_r}^{2k+1}e_{{2k,2k}}\big(2k
({r}+1)-r\big)F^{2k({r}+1)+{r}}_{2k({r}+1)+{r}}},
\\\label{k4}
{[\mathcal{E}^{0}_{2k+1,r}, {\mathbb{F}_{r}}]}&=&
\frac{{-{a_{r}}}^{2k+2}(2k+3)(2k+1)^{2k+1}_{-2}
({r}+1)}{(2k)!2^{2k+1}\big((2k+2)(r+1)+1\big)}E^{(2k+1)({r}+1)+{r}}_{(2k+1)({r}+1)+{r}}
\\\nonumber
&&-2{{a_{r}}}^{2k+2}{e}_{{2k+1,2k+2}}\big(2(k+1)({r}+1)+1\big)F^{{(2k+1)({r}+1)+{r}}}_{(2k+1)({r}+1)+{r}}.
\ea
On the other hand \cite[Lemma 3.4]{GazorMokhtariInt} implies that
\bes
\Theta^{(2k+1)({r}+1)+{r}}_{(2k+1)({r}+1)+{r}} \hbox{ and } F^{(2k+1)({r}+1)+{r}-1}_{(2k+1)({r}+1)+{r}-1}
\ees
are simplified in the \((r+1)\)th level of hypernormalization. For any sufficiently large number \(m,\) the sequence \(e_{2k+1,m}\) alternatively changes its sign and thus, \(e_{2k+1,2k+2}\) is nonzero. Thereby,
\bes
\hbox{ either  } E^{{(2k+1)({r}+1)+r}}_{(2k+1)(r+1)+r} \hbox{ or } F^{{(2k+1)(r+1)+r}}_{(2k+1)(r+1)+r}
\ees
is eliminated from the \((r+1)\)th level normal form, depending on whether style I or style II is applied. Now we show that \(e_{2k,2k}\) is nonzero by an induction on \(m\) through which the sign of \(e_{2k,m}\) is discussed.
Therefore by Equation \eqref{k4}, \(F^{2k({r}+1)+{r}}_{2k({r}+1)+{r}}\) can be removed from system and the proof is complete.
Since \(e_{2k,1}=-\frac{r}{2(r+2k+2)}<0,\) we have \(e_{2k,m}<0\) for \(0\leq m\leq k+1.\) Then,
 \bas
e_{2k,k+2}&=&{\frac{-(r+2)e_{2k,k+1}}{2\big((k+2)(r+1)+1\big)}}, \hbox{ and } e_{2k,k+2}>0.
\eas Next, \(e_{2k, k+3}<0,\) \(e_{2k, k+4}>0,\) and so forth. This concludes that \(e_{2k,2k}\) is nonzero.
\epr

\begin{lem}\label{YEF}
For any natural numbers \(m\) and \(n,\) there exist \(\mathfrak{F}^m_{n}, \mathfrak{E}^m_{n}\) such that
\ba\label{414}
{[\mathfrak{{F}}^m_n,\mathbb{F}_r]}+F^m_n&=\!&\!{\frac{{(-a_r)}^{n-m}(r(n-m)+n-2r)\big((m-n)(r+1)+m+2\big)_{2(r+1)}^{n-m-1}}
{{2}^{n-m}(m+2)(m+r+3)^{n-m-1}_{r+1}}}F^{rn-rm+n}_{rn-rm+n}\nonumber,
\ea
\ba\label{415}
\nonumber
{{[\mathfrak{{E}}^m_n,\mathbb{F}_r]+E^m_n}}&=\!\!&\!
{a_r}^{n-m}\left(
\frac{(r)^2_2\varphi_{n-m-1}}{n+r+2+r(n-m-1)}+\big(r(m-n+2\big)-n)\phi_{n-m-1}\right)\,F^{rn-rm+n}_{rn-rm+n}\\
&\!&\!+\frac{{a_r}^{n-m}(n+2)\big((n-m-1)(r+1)-m-1\big)^{n-m}_{-2(r+1)}}{(m+1)2^{n-m}\big(r(n-m)+n+2\big)
(m+r+2)^{n-m-1}_{r+1}}\,E^{rn-rm+n}_{rn-rm+n},
\ea
where the sequence \(\phi_l\) follows \(\phi_0=\phi_{-1}=0,\)
\ba\nonumber
{{\phi}_{l+1}}&=&\frac{(r)^2_2}{2{\big(n+2+l(r+1)\big)}{{\big({m+2+(l+1 )(r+1)}\big)} }}\varphi_l
-\frac {\big((2l+m-n)(r+1)+m+2\big)}{2\big({m+2+( l+1 )(r+1)}\big)}\phi_l,
\ea
and
\ba
{\varphi}_l=\frac{(n+2)\big((n-m-1)(r+1)-m-1\big)^{l}_{-2(r+1)}}{2^{l+1}(m+1)( n+lr+2 )(m+r+2)^{l}_{r+1}}.
\ea
\end{lem}
\bpr Let
{\bas
{\mathfrak{{E}}^m_n}&:=&\sum _{l=0}^{n-m-1}{a_r}^l{\phi}_lF^{m+lr+l+1}_{n+lr}
+\sum _{l=0}^{n-m-1}{a_r}^l{\varphi}_lE^{m+lr+l+1}_{n+{lr}},\\
\mathfrak{{F}}^m_n&:=&\sum _{l=0}^{n-m-1}{\frac{{(-a_r)^l}
\big((m-n)(r+1)+m+2\big)_{2(r+1)}^{l}}{{2}^{l+1}(m+2)(m+r+3)^{l}_{r+1}}}F^{m+lr+l+1}_{n+{lr}}.
\eas}
Then, the proof is a straightforward calculation.
\epr
Define
\be r_1:= r, \ \ {r_2}:= \min\{i\,|\,a_i\neq 0, i>r\}, \ \  {p}_1:= \min\{j\,|\,c_j\neq 0,j\geq p\}\ee
where \(a_i\) and \(c_j\) are coefficients of the \((r+1)\)-th level. When \(s<r_2,\) we {use} the generator \(\mathcal{X}_{r}^{k+1}\) to remove \(F^{2k(r+1)+r+s}_{2k(r+1)+r+s}\)
from the \((r+1)\)-th level normal form. However, this merely represented a Lie symmetry for the system in \cite{GazorMokhtariInt}.

\begin{thm}
There exist invertible transformation sending (\ref{Eq1}) into the \((s+1)\)-th level normal form
\ba\label{SNFCase1}
\left\{
  \begin{array}{ll}
\frac{dx}{dt}={2(y^2+z^2)}+a_rx^{r+1}+\sum_{k=0}^{\infty}(\alpha_{k+r_1}x^{r_1}+\beta_{k+s}x^{s})x^{k+1},
&\\
\frac{dy}{dt}=z-\frac{a_r(r+1)}{2}x^ry-\frac{y}{2}\sum_{k=0}^{\infty}\big({\alpha_{k+r_1}(k+r_1+1)}x^{r_1}
-{\beta_{k+s}}x^{s}\big)x^k+z\sum_{k=0}^{\infty}\gamma_{k+p_1} x^k,
  &\\
  \frac{dz}{dt}=-y-\frac{a_r(r+1)}{2}x^rz-\frac{z}{2}\sum_{k=0}^{\infty}\big({\alpha_{k+r_1}(k+r_1+1)}x^{r_1}
-{\beta_{k+s}}x^{s}\big)x^k-y\sum_{k=0}^{\infty}\gamma_{k+p_1} x^k,
  \end{array}
\right.
\ea
where \(s<r_2,\) we have
\begin{itemize}
 \item for the case of style I, \( \alpha_{k+r_1}=0\) for \({k\equiv_{2r+2} {r}}\) and \({k\equiv_{2r+2}{0}}\),
 \({k\equiv_{2r+2} {s}}\), and
  \(\beta_{k+s}=0\)  for \({k\equiv_{2r+2} {-(s+1)}}\), and  \( \gamma_{k+p_1}=0\) for \({k\equiv_{2(r+1)}{-(p_1+1)}}.\)
  \item for style II,  \( \alpha_{k+r_1}=0\) for \({k\equiv_{2r+2}{-1}}\), \({k\equiv_{2r+2} {0}}\), \({k\equiv_{2r+2} {s}}\), and \({k\equiv_{2r+2} {-(r+1)}}\) and \( \gamma_{k+p_1}=0\) for \({k\equiv_{2(r+1)}{-(p_1+1)}}.\)
\end{itemize}
When \(s\geq r_2\) the following holds:
\begin{itemize}
 \item style I: \( \alpha_{k+r_1}=0\) for \({k\equiv_{2r+2} {r}}\) and \({k\equiv_{2r+2}{0}}\), and \(\beta_{k+s}=0\)  for \({k\equiv_{2r+2} {-(s+1)}}\), where  \( \gamma_{k+p_1}=0\) for \({k\equiv_{r+1}{-(p_1+1)}}.\)
  \item style II:  \( \alpha_{k+r_1}=0\) for \({k\equiv_{2r+2}{-1}},\) \({k\equiv_{2r+2} {0}}\), and \({k\equiv_{2r+2} {-(r+1)}}\), while \( \gamma_{k+p_1}=0\) for \({k\equiv_{2(r+1)}{-(p_1+1)}}.\)
\end{itemize}
Besides, Equations \eqref{SNFCase1} represent an infinite level normal form.
\end{thm}
\bpr
Assume that \(s<r_2\) and the grading function follows Equation \eqref{grad11}. By Theorem \ref{YEF} there exist state solutions \(\mathfrak{E}^m_n, \mathfrak{F}^m_n\) for \((m,n)=\big(j(r+1)+s-1,2k+rj+s\big)\) such that
\bas
&&{[\mathfrak{E}^m_n+\mathfrak{F}^m_n,\mathbb{F}_r]+[\mathcal{X}^{k+1}_r, E^s_s]}
\\&=&\Bigg(\sum_{j=0}^{k+1}{{{k+1}\choose{j}}} {a_r}^{2k+1}
\bigg(
{\frac{(-1)^{j+1}(2k+rj)^2_2\big(2(j-k)(r+1)+s-r\big)_{2(r+1)}^{2k-j+1}}
{2^{2k-j+1}(2k+rj+s+2)\big(j(r+1)+s+1\big)^{2k-j}_{r+1}}}
\\&&+ {\frac {{\varphi}_{2k-j}(s)^2_2(r)^2_2(2k-j+2)}{(2k(r+1)+r+s+2)(2k+rj+s+2)}}\bigg)-\frac{{\phi}_{{2k-j}}\big(2k(r+1)+s-r\big)}{(2k+rj+s+2)} \Bigg)F^{2k(r+1)+r+s}_{2k(r+1)+r+s}.
\eas
So, \(F^{2k(r+1)+r+s}_{2k(r+1)+r+s}\) can be eliminated from the \((s+1)\)-th level of normalization for any \(k\in \N\).
For \(s\geq r_2,\) the transformation generators {\(\mathcal{X}^{k+1}_r\)} and \(\mathcal{T}^0_{2k,r}\) are extended to a symmetry of the system and therefore, they can not simplify the system any further; see \cite[Theorem 3.5]{GazorMokhtariInt}.
\epr
\subsection{Case ii: \(s<r\).}\label{caseii}

For this case we assume that
\bas
v^{(2)}:= F^{-1}_0+{b_s}{E^s_s}+\sum_{k=r}^{\infty}a_k{F^k_k}+
\sum^\infty_{k=s+1}b_k{E^k_k}+\sum^{\infty}_{k=1} c_{k}\Theta_k^k.
\eas
The grading function is defined by
\be\label{grad12}\delta(F^l_k)=
\delta(E^l_k)=s(k-l)+k\, \hbox{  and } \delta(\Theta^l_k)=s(k-l)+k+s+1.\ee Denote
\({\mathbb{E}_s}:=F^{-1}_0+b_{s}E^{s}_{s}.\)
\begin{lem}
The \((s+1)\)-th level normal form of (\ref{Eq1}) is
\be v^{(s+1)}= F^{-1}_0+{{b_s}}{E^s_s}+{a_r}^{(s+1)}{F^r_r}+\sum_{k>r} a_k^{(s+1)}{F^k_k}+
\sum_{k>s} b_k^{(s+1)}{E^k_k}+\sum_{k\geq1}{c}^{(s+1)}_k \Theta^k_k, \ee
where in
\begin{itemize}
  \item style I,\,\,\( b_k^{(s+1)}=0\) for \({k\equiv_{s+1} {2s}}\) and \({k\equiv_{s+1}  s}\) and \(k\neq s^2+2s\), while \( c_k^{(s+1)}=0\) for \({k\equiv_{ (s+1)} { s}} .\)
  \item style II,\,\,\( a_k^{(s+1)}=0\) when \(k\equiv_{s+1}  2s,\)
   \( b_k^{(s+1)}=0\)  for \(k\equiv_{s+1} s\) and \(k\neq s^2+2s\), while \( c_k^{(s+1)}=0\) if \(k\equiv_{s+1}  s.\)
\end{itemize}
\end{lem}
\bpr
Let  \be\label{gamma2}\mathscr{G}:={\ad_{(F^{1}_{0})}\circ \ad_{({\mathbb{E}_s})}}.\ee Then,
\bas
\mathscr{Gamma}(F^{l}_{k})& =&-4(l+1)(k-l+2) F^{l}_{k}
-{\frac{2{b_s}k(k+2)(k-l+1)}{k+s+2}} F^{s+l+1}_{k+s}\\&&+{\frac{2{b_s}s(s+2)(k-l)(k-l+1)}{k+s+2}}E^{s+l+1}_{k+s},
\\
\mathscr{Gamma} (E^{l}_{k})&=&-4 l(k-l+1)E^{l}_{k}-2{b_s}(k-s)(k-l)E^{s+l+1}_{k+s},
\\
\mathscr{Gamma}(\Theta^{l}_{k})&=&-4l(k-l+1)\Theta^{l}_{k}-2b_{s}k (k-l)
\Theta^{s+l+1}_{k+s}. \eas Then, \(\ker (\mathscr{G})= \Span\{\mathcal{F}^{-1}_{k}, {\mathcal{T}}^0_{k},  \mathcal{E}^0_{k}\}\) where
\ba\label{p11}
\mathcal{F}^{-1}_{k}&:=&{\sum_{m=0}^{k+1}\frac{(-b_{s})^{m}(k+2)(k)^m_s}
{{(2s+2 )}^{m}(m)!(ms+k+2)}F^{m(s+1)-1}_{ms+k}
-}\sum_{{m=0}}^{k+1}{b_{s}}^{m}h_{m}E^{m(s+1)-1}_{ms+k},\\
 \mathcal{E}^0_{k}&:=&\sum _{m=0}^{k}\frac{{(-b_{s})}^m(k-s)^m_s}{(2s+2)^mm!}E^{m(s+1)}_{ms+k},\,\qquad\qquad\qquad\qquad\qquad\qquad\,\,
\\\label{p112}
{\mathcal{T}}^0_{k}&:=&\sum _{m=0}^{k}{\frac{(-b_{{s}})^m(k)^m_s}{{(2s+2)^mm!}}}\Theta^{m ( s+1 )
}_{k+ ms}
,\,\,\qquad\qquad\qquad\qquad\qquad\qquad\qquad\ea
for \(h_0=0\) and
\bas
h_{m+1}  &=& {\frac{(-1)^{m+1}(s)^2_2(k)^{m}_s(k+2)(k-m+2)}{2^{m+1}{(s+1)}^{m}m!(ms+k+2)_s^2\big((m+1)(s+1)-1\big)}}
{-\frac{s(m-1)+k}{2\big((m+1)(s+1)-1\big)}h_{m}.
}\eas
Note that we have \(\delta(\mathcal{F}^{-1}_{{k}})\equiv_{s+1} {-1},\) \(\delta(\mathcal{E}^{0}_{{k}})\equiv_{s+1} 0,\) \(\delta( {\mathcal{T}}^0_{k})\equiv_{s+1} 0,\)
and \(\mathcal{E}^0_{s}= E^0_s.\) On the other hand
\ba\nonumber
{[\mathcal{F}^{-1}_{k},{\mathbb{E}_s}]}&= &{\frac{(-b_{s})^{k+2}\big(k(s+1)+s\big)_2^2{(k)}^{k+1}_s}{2^{k+1}{(s+1 )}^{k+2}(k+1)!\big((k+1)(s+1)+1}\big)F^{ks+k+2s}_{ks+k+2s}}
\\\label{Lp11}&& -{ 2({b_s})^{k+2}h_{k+2} \big((k+2)(s+1)-1\big)  E^{ks+k+2s}_{ks+k+2s}},\\\label{Lp12}
 {[\mathcal{E}^0_{k},{\mathbb{E}_s}]}&=&\frac{{(-b_s)}^{k+1}
 {(k-s)^{k+1}_s}}{2^kk!(s+1)^k}E^{ks+k+s}_{ks+k+s},\\\label{Lp13}
 [{\mathcal{T}}^0_{k},{\mathbb{E}_s}]&=&
{ \frac{{(-b_s)}^{k+1}(k)^k_s}{2^k(s+1)^{k-1}(k-1)!}\Theta^{ks+k+s}_{ks+k+s}}.
\ea
Notice that \(h_{1}=\frac{-(s+2)(k+2)}{2(k+s+2)}<0.\) By an induction argument we may conclude that for any odd number \(k,\) the number \(h_k\) is negative and \(h_k\) is positive for any even number \(k\). Therefore, \(h_k\neq0\) for any \(k\) and by Equations (\ref{Lp11}--\ref{Lp13}), the proof is complete.
\epr

\begin{lem}\label{Trans}
For any natural numbers \(m\) and \(n\), there exist \(\mathfrak{A}^m_n\) and \(\mathfrak{B}^m_n\) such that
\bas
F^m_n+[\mathfrak{A}^m_n,\mathbb{E}_s]&=& \left({\frac {b_{{s}}{\zeta}_{{n-m-1}}s(s+2)}{n+sn-sm+2}}+b_{{s}}{\xi}_
{{n-m-1}}\big(s(m-n+2)-n\big)\right)
 E^{sn-sm+n}_{sn-sm+n}
 \\&&
+{\frac{(-{b_s})^{n-m}(n)^{n-m}_s{( n+2 )}}{2^{n-m}(m+2)^{n-m+1}_{s+1}}F^{sn-sm+n}_{sn-sm+n}},
\\
{E^m_n+[\mathfrak{B}^m_n,\mathbb{E}_s]}&=&
{\frac{{{b_s}}^{n-m}(2s-n-sn+sm)(s-n)^{n-m-1}_{-s} }{2^{n-m}(m+1)(m+s+2)^{n-m-1}_{s+1}}E^{sn-sm+n}_{sn-sm+n}},
\eas
where
\bas
\xi_l&:=&\frac{-(s+2){b_s}^l(s-n)_{-s}^{l}\left((n+ls+2)(m+2)^{l}_{s+1}-
(n+2)(m+s+2)^{l}_{s+1}\right)}{2^{l+1}(s-n)(n+ls+2)(m+s+2)_{s+1}^{l}(m+2)_{s+1}^{l}},
\\
\zeta_l&:=&\frac{(-{b_s})^l(n+2)(n)^{l}_{s}}{2^{l+1}(n+ls+2)(m+2)^{l+1}_{s+1}}.
\eas
 \end{lem}
 \bpr
The proof is straightforward by
\bas
{\mathfrak{B}^m_n}&:=&{\sum _{l=0}^{n-m-1}\frac{{b_s}^l(s-n)^l_{-s} }{2^{l+1}(m+1)(m+s+2)^l_{s+1}}E^{m+ls+l+1}_{n+ls}},
\\
\mathfrak{A}^m_n &:=&{\sum_{l=0}^{n-m-1} \xi_l}{E^{m+ls+l+1}_ {n+ls}} +{\sum_{l=0}^{n-m-1}\zeta_lF^{m+ls+l+1}_{ n+ls}
}.\eas
\epr
We only use style II in the following theorem. Assume that the \((s+1)\)-th level coefficients \(a_i\) and \(b_j\) are nonzero for some \(i, j\geq 1.\) Define
\bes s_1:=s, \ \ s_2:=\min\big\{j\,|\,b_j\neq 0 \hbox{ for }  j> s\big\}, \ \ r:=\min\big\{i\,|\, a_i\neq0, i\geq 1\big\}.
\ees

\begin{thm}
There exist invertible transformations sending (\ref{Eq1}) into the \((r+s+2)\)-th level normal form
\ba\label{SNFCase2}
\frac{dx}{dt}&=&2(y^2+z^2)+a_rx^{r+1}+{b_s}x^{s+1}+x\sum_{k=1}^{\infty}(\alpha_{k+r}x^{r}+\beta_{k+s}x^{s})x^{k},
\\\nonumber
\frac{dy}{dt}&=&z-\frac{a_r(r+1)}{2}x^ry+\frac{1}{2}{b_s}x^{s}y-\frac{y}{2}\sum_{k=1}^{\infty}\big({\alpha_{k+r}(k+r+1)}x^{r}
-{\beta_{k+s}}x^{s}\big)x^k+z\sum_{k=1}^{\infty}\gamma_kx^k,
\\\nonumber
\frac{dz}{dt}&=&-y-\frac{a_r(r+1)}{2}x^rz+\frac{1}{2}{b_s}x^{s}z-\frac{z}{2}\sum_{k=1}^{\infty}\big({\alpha_{k+r}(k+r+1)}x^{r}
-{\beta_{k+s}}x^{s}\big)x^k-y\sum_{k=1}^{\infty}\gamma_kx^k.
\ea Here
\(\alpha_k=0\) for \(k\equiv_{s+1}  2s,\) \( \beta_k=0\)  for \(k\equiv_{s+1} s\) where \(k\neq s^2+2s\), \( \gamma_k=0\) for \(k\equiv_{s+1}  s,\) and
      \begin{description}
        \item[if \({s_2}<r\),] we have \( \beta_k=0\) for \(k= s^2+{s_2}+s,\)
        \item[when \({s_2}=r\),] \( \alpha_k=0\) for \(k= s^2+{s_2}+s,\)
        \item[for \({s_2}>r\),] \( \alpha_k=0\) for \(k= s^2+r+s.\)
      \end{description} Furthermore, the differential system \eqref{SNFCase2} is indeed the infinite level normal form.
\end{thm}
\bpr
\noindent For  ${s_2}<r$, by Lemma \ref{Trans} there exists a $\mathfrak{B}^{s_2}_{s+{s_2}}$ such that

\bas
[\mathfrak{B}^{s_2}_{s+{s_2}},\mathbb{E}_s] +[E^0_s,E^{{s_2}}_{{s_2}}]=\frac{(-{b_s})^{s}(s_2-s)_s^{s+1}}{2^{s}({s_2}+1)({s_2}+s+2)^{s-1}_{s+1}}E_{s^2+{s_2}+s}^{s^2+{s_2}+s}.
\eas
Hence, in both styles we can remove \(E_{s^2+{s_2}+s}^{s^2+{s_2}+s}\) from the system. For ${s_2}=r$, by  Lemma \ref{Trans} there exist  $\mathfrak{A}^{s_2}_{s+{s_2}}$  and $\mathfrak{B}
^{s_2}_{s+{s_2}}$ such that
\bas
&&{[\mathfrak{A}^{s_2}_{s+s_2}+\mathfrak{B}^{s_2}_{s+{s_2}},\mathbb{E}_s]+[E^0_s,a_{s_2}F^{s_2}_{s_2}+ b_{s_2}E^{s_2}_{s_2}]}
={\frac{{a_{s_2}(-b_s)}^s(s_2)^2_2( {s_2+s} )^{{s}}_s}{2^s(s_2+2)^{{s}+1}_{s+1}}F^{s^2+s_2+s }_{s^2+s_2+s }}
\\
&&+\Big(\frac{{{ b_{{s_{{2}}}}(-b_s)}}^{s}(s_2)^{s}_{s}( s_{{2}}-s)}{2^{s}(s_2+1)(s_2+s+2)^{s-1}_{s+1}} +{\frac {a_{{s_{{2}}}}{{(-b_s)}}^{s}(s)_2^2{(s_2)^{s}_{s} } }{2^s(s_{{2}}+s+2){(s_2+s+2)^{s-1}_{s+1}}}}
+{\frac {b_{{s}}a_{s_2}\zeta_{{s-1}}(s)^2_2(s_2)^2_2 }{(s_2+s+2)^2_{s^2}}}
\\
&&-\frac{b_{{s}}a_{s_2}\xi_{{s-1}}(s_2)^2_2({s}^{2}+s_2-s)}{(s_2+s+2)}\Big)E^{s^2+s_2+s }_{s^2+s_2+s}.
\eas Therefore, \(F^{s^2+{s_2}+s }_{s^2+{s_2}+s}\) can be eliminated in style II.
\noindent When ${s_2}>r$,{ by Lemma \ref{Trans} there exist  $\mathfrak{A}^r_{s+r}$ and  $\mathfrak{B}^r_{s+r}$ such }that
\bas
&&{[\mathfrak{A}^r_{s+r}+\mathfrak{B}^r_{s+r},\mathbb{E}_s]+[E^0_s,F^{r}_{r}]}=
\bigg(
{\frac {b_{{s}}\zeta_{{s-1}}(r)^2_2(s)^2_2}{(r+s+2)(r+s+2+{s}^{2})}}+\frac{b_{
{s}}(r)^2_2\xi_{{s-1}}( s-r-{s}^{2})}{(r+s+2)}
\\
&+&\frac{{{(-b_s)}}^{s}(s)^2_2(r)^{s}_{s} }{2^{s}(r+s+2)(r+s+2)^{s-1}_{s+1}}
\bigg)E^{s^2+r+s }_{s^2+r+s }+
\frac{{(-b_s)}^s(r)^2_2( {r+s} )^{{s}}_s}{2^s(r+2)^{{s}+1}_{s+1}}F^{s^2+r+s }_{s^2+r+s}.
\eas Since the  coefficients of \(F^{s^2+r+s }_{s^2+r+s}\) is nonzero, the proof is complete.
\epr

\subsection{Case iii: \(r=s\).}\label{caseiii}
Let
\ba\label{SCN3}
 {v}^{(2)}:=
F^{-1}_0+{a_s} F^s_s+{b_s} E^s_s+\sum_{k=s+1}^{\infty}(a_kF^k_k+b_kE^k_k)
+\sum_{k=1}^{\infty}c_k\Theta^k_k, \ea and assume that \(\frac{a_s}{b_s}\) is a non-algebraic number. This section is similar to the case III of Baider and Sanders \cite{baidersanders} which is the most difficult case of the three. We essentially use the idea of Wang {\rm et. al.} \cite{LiZhangWang14,KokubuWang,wlhj}; they assume that a certain ratio be non-algebraic. This prevents that the ratio be a root for any polynomial that appears in the computations and the results readily follows. When the corresponding fraction is an algebraic number, both problems (the case in this subsection and the third case of \cite{baidersanders}) remain unsolved. Denote
\be\mathbb{X}_s:=F^{-1}_0+{a_s} F^s_s+{b_s} E^s_s,\ee where  \(\delta(\mathbb{X}_s)=s\) and define a grading function by
\be\label{grad12}\delta(F^l_k):=
\delta(E^l_k):=s(k-l)+k\, \hbox{  and } \delta(\Theta^l_k):=s(k-l+1)+k+1.\ee

\begin{thm}\label{Thm4.7}
Let \(\frac{{a_s}}{{b_s}}\) be a non-algebraic number. Then, there exist invertible transformations that sends \(v^{(1)}\) given by equation (\ref{Eq1}) into the infinite level normal form system
\ba\label{SNFCase3}
\left\{
  \begin{array}{ll}
  \frac{dx}{dt}=2(y^2+z^2
)+({b_s}+{a_s})x^{s+1}+x\sum_{k=s+1}^{\infty}(\beta_k+\alpha_k)x^{k},
  &\\
  \frac{dy}{dt}=z+\,\big({b_s}-{{a_s}(s+1)}\big)\frac{x^{s}y}{2}+\,\frac{y}{2}\sum_{k=s+1}^{\infty}\big({\beta_k}-{\alpha_k(k+1)}\big)x^k
  +z\sum_{k=1}^{\infty}\gamma_kx^k,
  &\\
  \frac{dz}{dt}=-y+\big({b_s}-{{a_s}(s+1)}\big)\frac{x^{s}z}{2}+\frac{z}{2}\sum_{k=s+1}^{\infty}({\beta_k}-{\alpha_k(k+1)}){x^k}
  -y\sum_{k=1}^{\infty}\gamma_kx^k,
 \end{array}
\right.
\ea
where in
\begin{itemize}
\item  style I,\,\, \(\beta_k=0\) for \(k\equiv_{s+1}  s\) and  \(k\equiv_{_{s+1}} 2s\), and for \(k\equiv_{s+1}  -1\) we have
\(\gamma_k=0\).
\item style II,\,\, \(\alpha_k=0\) for \(k\equiv_{s+1} s \) and \(k\equiv_{s+1}  2s \), while \(\gamma_k=0\) for \(k\equiv_{s+1}  -1\).
\end{itemize}
\end{thm}
\bpr
Define
\be\label{gamma3}\mathscr{Gamma}:=\ad_{F^{1}_{0}}\circ \ad_{\mathbb{X}_s},\ee and
 \bas
\mathscr{Gamma}(F^{l}_{k})&=&-4(l+1)(k-l+2)
F^{l}_{k}+{\frac{2b_{{s}}s(s+2)(k-l)(k-l+1)}{k+s+2}} E^{l+s+1}_{k+s}
\\&&-2(k-l+1)\left({a_s}\big((l-k)(s+1)+l-s\big)
+{\frac{b_{{s}}k(k+2)}{k+s+2}}\right)F^{l+s+1}_{k+s},
\\
\mathscr{Gamma}(E^{l}_{k})&=&-4l(k-l+1)E^{l}_{k}+\frac
{2{a_s}s(s+2)(k-l+1)}{k+s+2}F^{l+s+1}_{k+s}
\\&&-2(k-l)\left(\frac{{a_s}(k+2)\big(l-(k-l)(s+1)\big)}{k+s+2}+{b_s}
(k-s)\right)E^{l+s+1}_{k+s},
\\
\mathscr{Gamma}(\Theta^{l}_{k})&=&-4l(k-l+1)\Theta^l_{k}-2(k-l)\Big(k{b_s}
+a_{s}\big(l({s}+2)-k({s}+1)\big)\Big)\Theta^{l+s+1}_{k+s}. \eas
Hence,

\be
\ker(\mathscr{Gamma})= \Span \{\mathcal{F}^{-1}_{k}, \mathcal{T}^0_{k},\mathcal{E}^0_{k}\,|\, k\in \N\},
\ee where
\bas
\mathcal{F}^{-1}_{k}&:=&\sum_{m=0}^{k+1}u^{-1}_{m}{{b_s}}^m F^{m(s+1)-1}_{ms+k}+\sum _{m=0}^{k+1}w^{-1}_{m}{{b_s}}^mE^{m(s+1)-1}_{ms+k},\qquad\,\,\,\
\delta(\mathcal{F}^{-1}_{k}) \equiv_{s+1} {(-1)},
\\
 {\mathcal{E}^0_{k}}&:=&\sum _{m=0}^{k}u^{0}_{m}{{b_s}}^mF^{m(s+1)} _{m s+k}+\sum_{m=0}^{k}w^{0}_{m}{{b_s}}^mE^{m(s+1)}_{ms+k},\qquad\quad\qquad
\delta(\mathcal{E}^{0}_{k})\equiv_{s+1} 0,
\\
\mathcal{T}^0_{k}&:=&\sum_{m=0}^{k}\frac{\big(k(a_ss+a_s-{b_s})\big)^{m}_{-\big(2a_s(s+1)+s{b_s}\big)}}
{m!(2( s+1 ))^{m}}\Theta^{m(s+1)}_{ms+k},\quad\ \ \ \ \ \ \
\delta(\mathcal{T}^{0}_{k})\equiv_{s+1} 0,
\eas
\(w^{-1}_0=u^0_0=0\) and \(u^{-1}_0=w^0_0=1.\) The sequences \(u^i_m\) and \(w^i_m\) for \(i=0, -1\) and \(m>0\) follow the recurrence relations
\ba\label{1}
u^{-1}_{m+1}&:=&
-\frac{(ms+k)(ms+k+2)}{2(m+1)(s+1)\big(s(m+1)+k+2\big)}u^{-1}_{m}
-\frac{a_s\big(2(m-1)-k\big)}{2{b_s}(m+1)}u^{-1}_{m}
\\\nonumber
&&
+\frac{{a_s}s(s+2)}{2{b_s}(m+1)(s+1)\big(s(m+1)+k+2\big)}w^{-1}_{m},
\\\label{2}
w^{-1}_{m+1}&:=&
\frac{{a_s}(ms+k+2)\big((k-2m+1)(s+1)+1\big)
}{2{{b_s}}\big(m(s+1)+s\big)\big(s(m+1)+k+2\big)}w^{-1}_{m}-\frac{\big(s(m-1)+k\big)}{2\big(m(s+1)+s\big)}w^{-1}_{m}
\\\nonumber &&
+\frac{{s(s+2)(k-m+2)}}{2\big(m(s+1)+s\big)\big(s(m+1)+k+2\big)}u^{-1}_{m},
\ea
\ba\label{3}
u^{0}_{m+1}&:=&
-\frac{(ms+k)(ms+k+2)}{2\big((m+1)(s+1)+1\big)\big(s(m+1)+k+2\big)}u^{0}_{m}
\\\nonumber
&&
-\frac{a_s\big((2m-k)(s+1)-s\big)}{2{b_s}\big((m+1)(s+1)+1\big)}u^{0}_{m}
+\frac{{a_s}s(s+2)}{2{b_s}\big((m+1)(s+1)+1\big)\big(s(m+1)+k+2\big)}w^{0}_{m},
\\\label{4}
w^{0}_{m+1}&:=&
\frac{{a_s}(k-2m)(ms+k+2)
}{2{{b_s}}(m+1)\big(s(m+1)+k+2\big)}w^{0}_{m}-\frac{s(m-1)+k}{2(m+1)(s+1)}w^{0}_{m}
\\\nonumber &&
+\frac{{s(s+2)(k-m+1)}}{2(m+1)(s+1)\big(s(m+1)+k+2\big)}u^{0}_{m}.
\ea
On the other hand
\bas
{[\mathcal{F}^{-1}_{k}, \mathbb{X}_s]}&=&{2{{b_s}}^{k+2}\Big((k+2)(s+1)u^{-1}_{{k+2}}F^{ks+k+2s}_{ks+k+2s}+
\big(s(k+2)+k+1\big)w^{-1}_{{k+2}}E^{ks+k+2s}_{ks+k+2s}\Big),}
\\\label{E0k}
{[\mathcal{E}^0_{k},\mathbb{X}_s]}&=&{2{{b_s}}^{k+1}\Big(\big((k+1)(s+1)+1\big)
u^0_{{k+1}}F^{ks+k+s}_{ks+k+s}+(k+1)(s+1)w^0_{{k+1}}E^{ks+k+s}_{ks+k+s}\Big)},
\\
{[\mathcal{T}^0_{k},\mathbb{X}_s]}&=&\frac{({a_s}+{b_s})\big(k(a_ss+a_s-{b_s})\big)
^{k}_{-\big(2a_s(s+1)+s{b_s}\big)}}{{2^k(k-1)!( s+1 )^{k-1}}}\Theta^{ks+k+s}_{ks+k+s}.
\eas
We claim that \(w^0_{{m+1}}\) and \(u^0_{m+1}\) are nonzero polynomials in terms of \(\frac{a_s}{b_s}\). For the case of \(u^{-1}_{{k+2}}\) and \(w^{-1}_{{k+2}},\) the results are similar. Due to the recurrence relations for \(k\neq s,\) the constant term of \(w^0_{{m+1}}\) is given by \(\frac{(-1)^m(k-s)^m_s}{2^mm!}.\) For \(k=s,\) the coefficient of \((\frac{a_s}{b_s})^{m+1}\) in the polynomial \(w^0_{{m+1}}\) is
\bes\frac{(k+2)(k)^{m}_{-2}}{2^mm!(k+km+2)}.\ees
On the other hand, the coefficient of \((\frac{a_s}{b_s})\) in \(u^0_{m+1}\) is governed by \bes\frac{(-1)^{m+1}m!{k}^{m}(k+2)_k^m}{2^m(k+2)_{k+1}^m(2k+2)_{k}^m}.\ees Since \(\frac{a_s}{b_s}\) is assumed a non-algebraic number, \(w^0_{k+1},\) \(u^0_{k+1},\) \(u^{-1}_{{k+2}}\) and \(w^{-1}_{{k+2}}\) in Equation \eqref{E0k} are nonzero. Hence, the proof is complete by the rules established in styles I and II.
\epr
\begin{rem}\label{Nonalgebraic}
The condition \(\frac{{a_s}}{{b_s}}\) being  non-algebraic was an essential assumption in the proof of Theorem \ref{Thm4.7}. However, we usually truncate the normal form system up to certain degree, say \(N\). Then, we may only need that \(\frac{{a_s}}{{b_s}}\) be distanced from the zeros of the polynomials generated by \(u^{-1}_{{k+2}},\) \(w^{-1}_{{k+2}},\) \(u^0_{{l+1}},\) and \(w^0_{{l+1}}\) for \bes k=0, 1, \ldots, \left\lfloor\frac{N-3}{s+1}\right\rfloor\hbox{ and } l=0, 1, \ldots, \left\lfloor\frac{N}{s+1}\right\rfloor-1;\ees for an instance of this see Example \ref{KS}.
\end{rem}

\section{Examples}\label{sec4}

In this section we provide the formulas for the first few coefficients of the simplest normal forms in terms of the coefficients of the original system. These are very useful for practical applications. Next, our results are applied on the R\"{o}ssler and Kuramoto--Sivashinsky equations to demonstrate the applicability of our results. Derivations of formulas follow exact (no numerical approximation) symbolic implementation of the results in Maple. The derived formulas up to the classical normal form coefficients are consistent with the results of Algaba {\rm et. al.} \cite{AlgabaHopfZ}.
We thank E. Gamero for sending us their classical normal form program and formulas of up to degree five for this comparison.

Consider a differential system
\be\label{4NF}
\left(\begin{array}{ccc}\dot{x} \\\dot{y} \\\dot{{z}} \\\end{array}\right)=\left(\begin{array}{ccc}0& 0& 0 \\0&0 & 1 \\0&-1 & 0
\\\end{array}\right)\left(\begin{array}{c}
x\\y\\{z}\end{array}
\right)+\sum_{2\leq i+j+k}\left(\begin{array}{ccc}{a}_{ijk}\\b_{ijk}
\\c_{ijk}\end{array}\right)x^iy^jz^k,
\ee
where \(a_{{0,2,0}}+a_{{0,0,2}}\neq0.\) Denote its cubic-truncated simplest normal form of Equation (\ref{4NF}) by
\ba\label{INF}
\frac{dx}{dt}&=& 2\alpha_0\rho^2+(\alpha_1+\beta_1)x^2+(\alpha_2+\beta_2)x^3,
\\\nonumber
\frac{d\rho}{dt}&=& \frac{1}{2}(\beta_1-\alpha_1)x\rho+\frac{1}{2}(\beta_2-3\alpha_2)x^2\rho,
\\\nonumber
\frac{d\theta}{dt}&=& 1+\gamma_1x+\gamma_2x^4.
\ea
where
\bas
\alpha_0&:=&\frac{1}{4}(c_{{0,2,0}}+c_{{2,0,0}}),\\
\alpha_1&:=&\frac{1}{3}(c_{{0,0,2}}-a_{{1,0,1}}-b_{{0,1,1}}),\\
\beta_1&:=&\frac{1}{3}(2c_{{0,0,2}}+a_{{1,0,1}}+b_{{0,1,1}}),\\
\gamma_1&:=&\frac{1}{2}(b_{{1,0,1}}-a_{{0,1,1}}),\\
\eas \(\alpha_2, \beta_2\) and \(\gamma_2\) are given in the appendix.

\begin{exm}
Consider the R\"{o}ssler equation
\ba\nonumber
\dot{x}&=&bz-cx+xz,\\\label{Ros}
\dot{y}&=&z+ay,\\\nonumber
{\dot{z}}&=&-y-x.
\ea

\noindent This system has Hopf-zero singularity for two sets of parameters. One is for parameter values \be a=c, \ b=1, \ 0<a^2<2,\ee while the system associated with the other set of parameter values has a simple dynamics; see \cite{AlgabaHopfZ}. We denote the vector field corresponding to \eqref{Ros} by \(R_a\). Therefore, the classical normal form of Equations \eqref{Ros} is given by
\bas
R^{(1)}_a&=&
\Theta^0_0-{\frac {a}{2( {2-{a}^{2}})^{\frac{3}{2}}  }}{ F^{-1}_0}-{\frac {a (8{a}^{2}+11 )}{ 12( 2-{a}^{2} )^{\frac{5}{2}} }} { F^0_1}-
{\frac {a ( {a}^{2}+1 ) }{3( {2-{a}^{2}})^{\frac{3}{2}}}}{F^1_1}
-{\frac {(5{a}^{2}+2)a}{4(2-{a}^{2})^{\frac{5}{2}}}}{F^2_2}
\\&&-\frac{a}{3(2-{a}^{2})^{\frac{3}{2}}}E^0_1-{\frac {a}{3({2-{a}^{2}})^{\frac{1}{2}}}} E^1_1 -{\frac{a}{4(2-{a}^{2})^{\frac{3}{2}} }}E^2_2
-{\frac{({a}^{2}+1)}{2(2-{a}^{2})}}\Theta^1_1.
\eas
By computing the second level normal form, we notice that this system is among the case iii (see Subsection \ref{caseiii}). Then,
\bas
R^{(\infty)}_{\mp1}&=&\Theta^0_0\pm\frac{1}{2}{ F^{-1}_0}\pm\frac{2}{3}{ F^1_1}\pm\frac{55}{32}{ F^2_2}\pm\frac{1}{3} E^1_1\pm\frac{9}{32} E^2_2-\Theta^1_1,
\eas
and for \(a\neq\pm1\) the infinite level normal form is given by
\ba\label{INFRossler}
\frac{dx}{dt}&=& {\frac {-a{\rho}^{2}}{{({2-{a}^{2}})}^\frac{3}{2} }}-{
\frac {a{x}^{2}}{( 2-{a}^{2})^{\frac{3}{2}}}}-{\frac {a
 ( {a}^{2}+1 ) {x}^{3}}{ ( 2-{a}^{2} )^{\frac{5}{2}}}},
\\\nonumber
\frac{d\rho}{dt}&=& {\frac {{a}^{3}x\rho}{2( 2-{a}^{2})^{\frac{3}{2}}}}+{
\frac {3a(11{a}^{2}+2){x}^{2}\rho }{16( 2-{a}^{2})^{\frac{5}{2}}}},
\\\nonumber
\frac{d\theta}{dt}&=& 1-{\frac{({a}^{2}+1)}{2(2-{a}^{2})}}x.
\ea

\end{exm}

\begin{exm}\label{KS}
The standing waves of the Kuramoto--Sivashinsky equation gives rise to
\bas
\dot{x}&=&y,\\
\dot{y}&=&z,\\
\dot{z}&=&\mu x-2x^2-y,
\eas which has a Hopf-zero singularity at the origin; see Chang \cite{ChangKS}. Our formulas for second level normal form of up to degree six is given by
\bas
K^{(2)}=\Theta^0_0-\frac{1}{2}F^{-1}_0-2F^1_1
+{\frac {19}{3}}\Theta^2_2+{\frac {343}{9}}F^3_3-{\frac {6752297}{3240}}F^5_5-{
\frac {57385}{108}}\Theta^4_4.
\eas
This falls within the case i (see Subsection \ref{casei}). Hence, the infinite level normal form truncated up to degree six is governed by
\ba
K^{(\infty)}={\Theta^0_0-\frac{1}{2}F^{-1}_0-2F^1_1+\frac {686}{3}}E^3_3+{\frac {19}{3}}\Theta^2_2
-{\frac {57385}{108}}\Theta^4_4,
\ea or equivalently,
\ba\nonumber
\dot{x}&=&-{\rho}^{2}-2{x}^{2}+{\frac {686}{3}}{x}^{4},\\
\dot{\rho}&=&2x\rho+{\frac {343}{3}}{x}^{3}\rho,\\\nonumber
\dot{\theta}&=& 1+{\frac {19}{3}}{x}^{2}-{\frac {57385}{108}}{x}^{4}.
\ea
\end{exm}

\section{ Appendix }

In this appendix we present the formulas derived by symbolic implementation of the results in Maple. Denote \(\alpha_+:=\alpha_2\) and \(\alpha_-:=\beta_2.\) Then,
\bas
\alpha_{\pm}&=&\frac{-1}{4(a_{{0,2,0}}+a_{{0,0,2}})}\bigg((a_{{0,2,0}}+a_{{0,0,2}})
\Big(c_{{2,0,1}}+b_{{2,1,0}}-a_{{3,0,0}}-3\,c_{{2,0,0}}a_{{1,1,0}}+3\,b_{{2,0,0}}a_{{1,0,
1}}
\\
&&+c_{{2,0,0}}c_{{0,1,1}}+2c_{{2,0,0}}b_{{0,2,0}}\Big)\pm\frac{1}{48}(2a_{{2,0,0}}+c_{{1,0,1}}+b_{{1,1,0}})
\Big(4b_{{0,1,1}}b_{{0,0,2}}+8a_{{1,2,0}}+8a_{{1,0,2}}
\\&&
-12c_{{0,0,3}}+6c_{{1,1,0}}a_{{0,2,0}}-6c_{{1,1,0}}a_{{0,0,2}}-4c_{{0,2,0}}c_{{0,1,1}}-8c_{{0,2,0}}b_{{0,2,0}}
+6c_{{1,0,1}}a_{{0,1,1}}-4b_{{0,1,2}}
\\&&
+8\,c_{{0,0,2}}b_{{0,0,2}}-6\,b_{{1,1,0}}a_{{0,1,1}}+4\,b_{{0,2,0}}b_{{0,1,1}}+6\,b_{{1,0,1}}
a_{{0,2,0}}-6\,b_{{1,0,1}}a_{{0,0,2}}-4\,c_{{0,1,1}}c_{{0,0,2}}
\\&&
-12b_{{0,3,0}}
+8\,c_{{0,2,0}}a_{{1,1,0}}+8\,c_{{0,0,2}}a_{{1,1,0}}-8\,b_{{0,2,0}}a_{{1,0,1}}-8\,b_{{0,0,2}}a_{{1,0,1}}\Big)
\mp\frac{1}{3}(a_{{2,0,0}}-c_{{1,0,1}}
\\&&
-b_{{1,1,0}})
\Big(
3c_{{0,0,3}}
-b_{{0,1,1}}b_{{0,0,2}}
+c_{{0,2,0}}c_{{0,1,1}}+{2}c_{{0,2,0}}
b_{{0,2,0}}
+c_{{0,1,1}}c_{{0,0,2}}
-2c_{{0,0,2}}
b_{{0,0,2}}+
b_{{0,1,2}}
\\&&
-b_{{0,2,0}}b_{{0,1,1}}
+
a_{{1,0,2}}
+b_{{0,3,0}}
+c_{{0,2,0}}a_{{1,1,0}}+c_{{0,0,2}}a_{{1,1,0}}-b_{{0,2,0}
}a_{{1,0,1}}-b_{{0,0,2}}a_{{1,0,1}}
+a_{{1,2,0}}
\Big)\bigg),
\eas and
\bas
\gamma_2&:=&-\frac{1}{2}c_{{2,1,0}}-c_{{2,0,0}}c_{{0,2,0}}+\frac{1}{2}c_{{2,0
,0}}b_{{0,1,1}}-c_{{2,0,0}}a_{{1,0,1}}-\frac{1}{8}{c_{{1,1,0}}}^{2}
-\frac{1}{4}c_{
{1,1,0}}b_{{1,0,1}}-\frac{1}{8}{c_{{1,0,1}}}^{2}
\\&&+\frac{1}{4}c_{{1,0,1}}b_{{1,1,0}}
+\frac{1}{2}c_{{0,1,1}}b_{{2,0,0}}-b_{{2,0,0}}b_{{0,0,2}}-b_{{2,0,0}}a_{{0,1
,1}}-\frac{1}{8}{b_{{1,1,0}}}^{2}-\frac{1}{8}{b_{{1,0,1}}}^{2}+\frac{1}{2}b_{{2,0,1}}.
\eas

\end{document}